\documentclass[12pt]{article}
\usepackage[noend]{algpseudocode}
\usepackage{algorithmicx,algorithm}
\usepackage{amsmath}
\usepackage{indentfirst}
\usepackage{exscale}
\usepackage{relsize}
\usepackage{fancyhdr}  
\usepackage{amssymb}
\usepackage{mathrsfs}
\usepackage{color}
\usepackage{cases}
\usepackage{amsfonts}
\usepackage{enumerate}
 \usepackage{eucal}
\usepackage{bbding}
\usepackage[numbers,sort&compress]{natbib}
\allowdisplaybreaks[3]
\newtheorem{theorem}{Theorem}[section]

\newtheorem{lemma}[theorem]{Lemma}

\numberwithin{equation}{section}
  \usepackage{a4wide}

\def \O{{\Omega}}

\def\O{{\mathcal O}}

\begin{document}
\baselineskip 14pt
\title{\bf Note on a sum involving the divisor function}
\author{{Liuying  Wu}\\{ School of Mathematics and Statistics,}\\{ Hunan University of Science and Technology}
\\ {Xiangtan, 411201, Hunan, P. R. China.}\\{Email: liuyingw@tongji.edu.cn}}

\date{}
\maketitle
\noindent {\bf Abstract}:  Let $d(n)$ be the divisor function and denote by $[t]$ the integral part of the real number $t$. In this paper, we prove that
$$\sum_{n\leq x^{1/c}}d\left(\left[\frac{x}{n^c}\right]\right)=d_cx^{1/c}+\O_{\varepsilon,c}
\left(x^{\max\{(2c+2)/(2c^2+5c+2),5/(5c+6)\}+\varepsilon}\right),$$
where $d_c=\sum_{k\geq1}d(k)\left(\frac{1}{k^{1/c}}-\frac{1}{(k+1)^{1/c}}\right)$ is a constant. This result constitutes an improvement upon that of Feng.


\section{\bf Introduction} \label{s1}

Let $[x]$ denote the integral part of $x$, and let $d(n)$ be the divisor function. From \cite{Bour}, we know that  
$$\sum_{n\leq x}\left[\frac{x}{n}\right]=x\log x+(2\gamma-1)x+\O(x^{517/1648+o(1)}),$$
where $\gamma=0.577215664\cdots$ is the Euler constant. Recently, there has been a great deal of interest in estimating sums of the form
$$S_f(x):=\sum_{n\leq x}f\left(\left[\frac{x}{n}\right]\right),$$
where $f$ is an arithmetic function. In \cite{Bord}, Bordellès, Dai, Heyman, Pan and Shparlinski proved that if there exists $\alpha\in (0,2)$ such that 
$$\sum_{n\leq x}|f(n)|^2\ll x^\alpha,$$
then
$$S_f(x)=x\sum_{n\geq 1}\frac{f(n)}{n(n+1)}+\O\left(x^{(\alpha+1)/3}(\log x)^{(\alpha+1)/3+o(1)}\right).$$
Later, this estimate was improved by Wu \cite{Wu}, who proved that if there is a constant $\vartheta\in[0,1)$ such that $|f(n)|\ll n^\vartheta$ for $n\geq1,$ then
\begin{equation}\label{e.101}
S_f(x)=x\sum_{n\geq 1}\frac{f(n)}{n(n+1)}+\O\left(x^{(\vartheta+1)/2}\right).
\end{equation}
In particular, applying \eqref{e.101} to the von Mangoldt function $\Lambda(n),$ we have
$$S_\Lambda(x)=x\sum_{n\geq 1}\frac{\Lambda(n)}{n(n+1)}+\O\left(x^{1/2+\varepsilon}\right).$$
The $1/2-$barrier was first broken by Ma and Wu. In \cite{MW}, they showed that
$$S_\Lambda(x)=x\sum_{n\geq 1}\frac{\Lambda(n)}{n(n+1)}+\O\left(x^{35/71+\varepsilon}\right).$$
The exponent $35/71$ was improved successively to $97/203,9/19,92/195$ by Bordellès \cite{Bord1}, Liu, Wu and Yang \cite{LWY} and Zhang \cite{Zhang} respectively.
\vskip 2 pt
On the other hand,  using the fact that $|d(n)|\ll n^\varepsilon$ for $n\geq1$, we can get from \eqref{e.101} that
\begin{equation*}
S_d(x)=x\sum_{n\geq 1}\frac{d(n)}{n(n+1)}+\O\left(x^{1/2+\varepsilon}\right).
\end{equation*}
By a similar method that used in \cite{MW}, Ma and Sun \cite{MS} proved that
$$S_d(x)=x\sum_{n\geq 1}\frac{d(n)}{n(n+1)}+\O\left(x^{11/23+\varepsilon}\right).$$
Later, the exponent $11/23$ was improved successively to $19/40,5/11$ by Bordellès \cite{Bord1} and Stucky \cite{Stu} respectively.
\vskip 2 pt
Very recently, Feng \cite{Feng} studied the asymptotic formula of the sum
\begin{equation}\label{e.102}
S_{d,c}(x)=\sum_{n\leq x^{1/c}}d\left(\left[\frac{x}{n^c}\right]\right),
\end{equation}
where $c$ is a positive real number. Feng obtained
$$S_{d,c}(x)=x^{1/c}\sum_{k\geq1}d(k)\left(\frac{1}{k^{1/c}}-\frac{1}{(k+1)^{1/c}}\right)
+\O_{\varepsilon,c}\left(x^{\theta_c+\varepsilon}\right)$$
for any $\varepsilon>0,$ where 
\begin{subnumcases}{\theta_c=}
    \frac{2}{3c+2}\ & \text{if}~~$0<c<\frac{2}{11}$,\nonumber\\
    \frac{11}{11c+12}\ & if~~$c\geq \frac{2}{11}$.\nonumber
\end{subnumcases}
\vskip 2 pt
Motivated by the key ideas used in \cite{Stu}, we shall continue to study \eqref{e.102} and prove the following result.
\begin{theorem}
For any $\varepsilon>0$, we have
$$S_{d,c}=x^{1/c}\sum_{k\geq1}d(k)\left(\frac{1}{k^{1/c}}-\frac{1}{(k+1)^{1/c}}\right)
+\O_{\varepsilon,c}\left(x^{\theta_c+\varepsilon}\right),$$
as $x\rightarrow\infty,$ where
\begin{subnumcases}{\theta_c=}
    \frac{2c+2}{2c^2+5c+2}\ & \text{if}~~$0<c<\frac{2}{3}$,\nonumber\\
    \frac{5}{5c+6}\ & if~~$c\geq \frac{2}{3}$.\nonumber
\end{subnumcases}
\end{theorem}
\vskip 2 pt

\noindent{\bf Remark.} In order to compare our result with previous works, we note that 
$$\frac{11}{11c+12}>\max\left\{\frac{2c+2}{2c^2+5c+2},\frac{5}{5c+6}\right\}$$
for $c>2/9,$ which implies Theorem 1.1 improves Feng's error term if $c>2/9.$ It is also worth mentioning
that $\theta_c=5/11$ when $c=1,$ thus our result generalize the work of \cite{Stu}.
\vskip 2 pt
{\bf Notation.} Throughout this paper, we denote by $[x]$ for the integral part of any real number $x$, and let $\psi(t):=x-[x]-1/2.$ We use $\varepsilon$ to denote a sufficiently small positive number, and the value of $\varepsilon$ may change from statement to statement. We write $f=\mathcal{O}(g)$ or, equivalently, $f \ll g$ if $|f|\leq Cg$ for some positive number $C$. And we use $e(t)$ to denote $e^{2\pi it}.$

\section{Two important lemmas}
In order to prove Theorem 1.1, we need the following two important lemmas, the first one is due to Vaaler \cite[Theorem A.6]{GK}, and the second is from Jutila \cite[Lemma 4.6]{Jut}.

\begin{lemma}
For $x\geq1$ and $H\geq1$, we have
$$\psi(x)=-\sum_{1\leq |h|\leq H}\Phi\left(\frac{h}{H+1}\right)\frac{e(hx)}{2\pi ih}+R_{H}(x),$$
where $\Phi(t)=\pi t(1-|t|)\cos(\pi t)+|t|$, and
$$|R_{H}(x)|\leq\frac{1}{2H+2}\sum_{|h|\leq H}\left(1-\frac{|h|}{H+1}\right)e(hx).$$
\end{lemma}

\begin{lemma}
Let $2 \leqslant M < M ^ { \prime } \leqslant 2 M $, and let $f$ be a holomorphic function in the domain
$$D = \{ z : |z-x|<cM~for~some~x \in \left[M,M'\right]\},$$
where $c$ is a positive constant. Suppose that $f(x)$ is real for $M \leqslant x \leqslant M ^ { \prime } $, and that either 
\[f ( z ) = B z ^ { \alpha } \left( 1 + \O\left( F^{-1/3} \right) \right) \quad for~z \in D ,\] 
where $\alpha \neq 0 , 1$ is a fixed real number, and 
\[F =|B|M ^ { \alpha } , \quad o r \quad f ( z ) = B \log z \left( 1 + o\left( F ^ { - 1 / 3 } \right) \right) \quad for~z \in D ,\] 
where $F = |B| $.
 
Let $g \in C^{1}\left[M, M^{\prime}\right] ,$ and suppose that $M \leqslant x \leqslant M ^ { \prime } $, 
\[| g ( x ) | \ll G , \quad \left|g^{\prime}(x)\right| \ll G^{\prime}.\] 
Suppose also that $M^{3/4} \ll F\ll M^{3/2}$, then 
\[\left| \sum _ { M \leqslant m \leqslant M^{ \prime } } d ( m ) g ( m ) e ( f ( m ) ) \right| \ll \left( G + FG^{\prime} \right) M^{1/2} F^{ 1 / 3 + \varepsilon }.\] 
\end{lemma}

\section{Proof of Theorem 1.1}

Let $N\in[1,x^{1/(c+1)})$ be a parameter to be chosen later. We split the sum $S_{d,c}(x)$ into two parts:
\begin{equation}\label{e.301}
S_{d,c}(x)=S_1(x)+S_2(x),
\end{equation}
where
$$S_1(x):=\sum_{n\leq N}d\left(\left[\frac{x}{n^c}\right]\right),\qquad S_2(x):=\sum_{N<n\leq x^{1/c}}d\left(\left[\frac{x}{n^c}\right]\right).$$
Using that $d(n)\ll_\varepsilon n^\varepsilon$ for all $n\geq1$ and any $\varepsilon>0$, we have trivially
\begin{equation}\label{e.302}
S_1(x)\ll_\varepsilon Nx^\varepsilon.
\end{equation}
\vskip 2 pt
In order to bound $S_2(x),$ we put $k=[x/n^c]$. Then
$$k\leq \frac{x}{n^c}<k+1\Leftrightarrow \left(\frac{x}{k+1}\right)^{1/c}<n\leq \left(\frac{x}{k}\right)^{1/c}.$$
Thus
\begin{align}\label{e.303}
S_2(x)&=\sum_{k\leq x/N^c}d(k)\sum_{\left(\frac{x}{k+1}\right)^{1/c}<n\leq \left(\frac{x}{k}\right)^{1/c}}1\nonumber\\
&=\sum_{k\leq x/N^c}d(k)\left\{\frac{x^{1/c}}{k^{1/c}}-\psi\left(\frac{x^{1/c}}{k^{1/c}}\right)
-\frac{x^{1/c}}{(k+1)^{1/c}}+\psi\left(\frac{x^{1/c}}{(k+1)^{1/c}}\right)\right\}\nonumber\\
&=x^{1/c}\sum_{k\geq1}d(k)\left(\frac{1}{k^{1/c}}-\frac{1}{(k+1)^{1/c}}\right)
+R_0(x)-R_1(x)+\O_\varepsilon(Nx^\varepsilon),
\end{align}
where the bounds
$$x^{1/c}\sum_{k>x/N^c}d(k)\left(\frac{1}{k^{1/c}}-\frac{1}{(k+1)^{1/c}}\right)\ll Nx^\varepsilon,$$
$$\sum_{k\leq N}d(k)\left\{\psi\left(\frac{x^{1/c}}{(k+1)^{1/c}}\right)-\psi\left(\frac{x^{1/c}}{k^{1/c}}\right)\right\}\ll Nx^\varepsilon$$
are used, and $R_\delta(x)$ is defined by
$$R_\delta(x):=\sum_{N<k\leq x/N^c}d(k)\psi\left(\frac{x^{1/c}}{(k+\delta)^{1/c}}\right)$$
with $\delta=0,1.$
\vskip 2 pt
Now it remains to estimate 
\begin{equation}\label{e.304}
R_\delta(x)\ll x^\varepsilon \max_{N<D\leq x/N^c}\mathfrak{S}_{\delta,c}(x,D),
\end{equation}
where
$$\mathfrak{S}_{\delta,c}(x,D)=\sum_{D<k\leq 2D}d(k)\psi\left(\frac{hx^{1/c}}{(k+\delta)^{1/c}}\right).$$
By Lemma 2.1 and noticing the fact that $0<\Phi(t)<1~(0<|t|<1)$, we can derive that
\begin{equation}\label{e.201}
\mathfrak{S}_{\delta,c}(x,D)\ll \frac{D}{H}+\sum_{h\leq H}\frac{1}{h}\left|\sum_{D<k\leq 2D}d(k)e\left(\frac{hx^{1/c}}{(k+\delta)^{1/c}}\right)\right|
\end{equation}
for $N<D\leq x/N^c$ and $1\leq H\leq D.$ 
\vskip 2 pt
Then we will focus on the estimate of the inner sum over $k$. By Lemma 2.2 with
$$M=D,\quad M'=2D,\quad f(z)=\frac{hx^{1/c}}{(z+\delta)^{1/c}},\quad g(x)=1,\quad B=hx^{1/c},\quad F=\frac{hx^{1/c}}{D^{1/c}},$$
we find that
\begin{equation}\label{e.202}
\sum_{D<k\leq 2D}d(k)e\left(\frac{hx^{1/c}}{(k+\delta)^{1/c}}\right)\ll D^{1/2-1/(3c)}h^{1/3}x^{1/(3c)+\varepsilon},
\end{equation}
where
\begin{equation}\label{e.203}
\frac{D^{3/4+1/c}}{x^{1/c}}\ll 1\leq h\leq H\ll \frac{D^{3/2+{1/c}}}{x^{1/c}}.
\end{equation}
By \eqref{e.201} and  \eqref{e.202}, we have
\begin{align}\label{e.204}
\mathfrak{S}_{\delta,c}(x,D)\ll DH^{-1}+D^{1/2-1/(3c)}H^{1/3}x^{1/(3c)+\varepsilon}
\end{align}
for all $H$ satisfying \eqref{e.203}. Inserting \eqref{e.204} into \eqref{e.304}, we get
\begin{equation}\label{e.305}
R_\delta(x)\ll x^\varepsilon \max_{N<D\leq x/N^c}\left(DH^{-1}+D^{1/2-1/(3c)}H^{1/3}x^{1/(3c)+\varepsilon}\right).
\end{equation}
\vskip 2 pt
For $0<c<\frac{2}{3}$, from \eqref{e.305} we deduce that
$$R_\delta(x)\ll x^\varepsilon \left(xH^{-1}N^{-c}+N^{1/2-1/(3c)}H^{1/3}x^{1/(3c)+\varepsilon}\right).$$
Inserting this estimate into \eqref{e.303}, and combining \eqref{e.302} and \eqref{e.303} with \eqref{e.301}, it follows that 
\begin{equation}\label{e.308}
S_{d,c}(x)=x^{1/c}\sum_{k\geq1}d(k)\left(\frac{1}{k^{1/c}}-\frac{1}{(k+1)^{1/c}}\right)+R(x),
\end{equation}
where
$$R(x)=\O_{\varepsilon,c}\left(Nx^\varepsilon+x^{1+\varepsilon} H^{-1}N^{-c}+N^{1/2-1/(3c)}H^{1/3}x^{1/(3c)+\varepsilon}\right).$$
Taking
$$H=x^{3/4-1/(4c)}N^{1/(4c)-3c/4-3/8}\quad \text{and}\quad N=x^{2(1+c)/(2c^2+5c+2)},$$
it is easy to see that the condition \eqref{e.203} holds for all $D\in (N,x/N^c]$. Then we obtain
\begin{equation}\label{e.306}
S_{d,c}(x)=x^{1/c}\sum_{k\geq1}d(k)\left(\frac{1}{k^{1/c}}-\frac{1}{(k+1)^{1/c}}\right)
+\O_{\varepsilon,c}\left(x^{2(1+c)/(2c^2+5c+2)+\varepsilon}\right).
\end{equation}
\vskip 2 pt
Similarly, by \eqref{e.305} we have
$$R_\delta(x)\ll x^\varepsilon \left(xH^{-1}N^{-c}+N^{1/3-c/2}H^{1/3}x^{1/2+\varepsilon}\right).$$
for $c\geq \frac{2}{3}$. In this case, the error term in \eqref{e.308} becomes 
$$\widetilde{R}(x)=\O_{\varepsilon,c}\left(Nx^\varepsilon+x^{1+\varepsilon} H^{-1}N^{-c}+N^{1/3-c/2}H^{1/3}x^{1/2+\varepsilon}\right).$$
Taking
$$H=x^{3/8}N^{-3c/8-1/4}\quad \text{and}\quad N=x^{5/(5c+6)},$$
also, the condition \eqref{e.203} holds for all $D\in (N,x/N^c]$. Thus, we get
\begin{equation}\label{e.307}
S_{d,c}(x)=x^{1/c}\sum_{k\geq1}d(k)\left(\frac{1}{k^{1/c}}-\frac{1}{(k+1)^{1/c}}\right)
+\O_{\varepsilon,c}\left(x^{5/(5c+6)+\varepsilon}\right).
\end{equation}
Now Theorem 1.1 follows from \eqref{e.306} and \eqref{e.307}.

\bigskip

\noindent {\bf Acknowledgement.} We wish to thank the referee for a thorough reading of the paper. 


\end{document}